\documentclass[preprint,12pt]{elsarticle}

\usepackage{amssymb}
\usepackage{amsmath}
\usepackage{amsthm}
\usepackage{graphicx}

\usepackage{float}

\usepackage{booktabs,siunitx}
\usepackage{multirow}

\usepackage{fullpage}
\usepackage{bm}
\usepackage{lineno}

\usepackage[utf8x]{inputenc}
\usepackage{color}
\usepackage{algpseudocode, algorithm, algorithmicx}

\usepackage{setspace}
\newcommand{\rev}[1]{\textcolor{black}{#1}}

\usepackage{romannum}

\biboptions{square,numbers,sort&compress}

\begin{document}
\pagenumbering{arabic}

\begin{frontmatter}

\title{Inhomogenous Regularization with Limited and Indirect Data}
\author{Jihun Han\corref{cor1}}
\ead{jihun.han@dartmouth.edu}
\cortext[cor1]{Corresponding author.}
\author{Yoonsang Lee}
\ead{yoonsang.lee@dartmouth.edu}

\address{Department of Mathematics, Dartmouth College, \\27 N. Main Street, Hanover NH,  03755 USA}

\begin{abstract}
For an ill-posed inverse problem, particularly with incomplete and limited measurement data, regularization is an essential tool for stabilizing the inverse problem. Among various forms of regularization, the $\ell_p$ penalty term provides a suite of regularization of various characteristics depending on the value of $p$. When there are no explicit features to determine $p$, a spatially varying inhomogeneous $p$ can be incorporated to apply different regularization characteristics that change over the domain. 
This study proposes a strategy to design the exponent $p$ when the first and second derivatives of the true signal are not available, such as in the case of indirect and limited measurement data. The proposed method extracts statistical and patch-wise information using multiple reconstructions from a single measurement, which assists in classifying each patch to predefined features with corresponding $p$ values. We validate the robustness and effectiveness of the proposed approach through a suite of numerical tests in 1D and 2D, including a sea ice image recovery from partial Fourier measurement data. Numerical tests show that the exponent distribution is insensitive to the choice of multiple reconstructions.
\end{abstract}


\journal{Journal of Computational and Applied Mathematics}

\end{frontmatter}

\section{Introduction}
Inverse problems aim to identify an underlying signal from measurement data. The measurement data is typically limited and incomplete, and thus the inverse problem is often formulated as an ill-posed problem.
A regularization is an essential tool in stabilizing an ill-posed inverse problem. Regularization imposes an additional structure in the solution or restricts the solution to a particular space so that there is a unique solution to the problem. A general approach of regularization adds a penalty term to an objective function that measures the data misfit. Tikhonov, or $\ell_2$ regularization, adds an $\ell_2$ penalty term so that all components of the solution are equally regularized. In the Bayesian context, the Tikhonov regularization is equivalent to prior knowledge that the uncertainty of the solution follows a Gaussian distribution, which helps recover a smooth solution.

Another widely used class of regularization includes $\ell_1$ regularization to recover a sparse solution. $\ell_1$ regularization has been used as a convex relaxation of the $\ell_0$ penalty that imposes a constraint in the number of nonzero components of the solution. By avoiding the combinatorial complexity related to the $\ell_0$ penalty term, $\ell_1$ regularization is computationally efficient in recovering a sparse signal from a limited measurement data that is much smaller than the unknown variable's dimension \cite{CandesTao}. In particular, the total variation (TV) regularization \cite{rudin1992nonlinear}, which is related to an $\ell_1$ norm of the gradient,  shows robust performance in recovering jump-discontinuities or edges by incorporating the sparsity in the gradient domain.

When there are no explicit features, such as sparsity for $\ell_1$ or smooth variations for $\ell_2$, or when the features are mixed, none of the single regularizations provide robust performance. $\ell_2$ regularization smoothes out sharp edges or discontinuities. On the other hand, $\ell_1$ or TV gives rise to an unfavorable staircase effect \cite{dobson1996recovery} that yields piecewise constant recovery of smooth or oscillatory signals.
There are approaches to combine both the Tikhonov and TV regularizations \rev{\cite{liu2014adaptive, gholami2013balanced, asadi2019data}}, which utilize the trade-off characteristics of both regularizations.
The approach in \cite{liu2014adaptive} decomposes the signal into a piecewise constant component and a smooth component and penalizes two components with TV and Tikhonov regularizations, respectively. \rev{The approach proposed in \cite{gholami2013balanced} uses a weighted sum of both regularizations,  in which the balance is determined by the local behavior estimated from gradient information.} \rev{Moreover,  the $\ell_1$ could be applied to the specific space for sparse features of interest  \cite{asadi2019data}.  The regularized space for the true signal could be learned from the data \cite{lunz2018adversarial,heaton2020wasserstein},  which is expressed in a neural network formulation of  regularization \cite{lunz2018adversarial},  or the learned projection operator is directly incorporated into an optimization method \cite{heaton2020wasserstein}.  }

\rev{Other methods to handle mixed features include the weighted TV \cite{candes2008enhancing,  chartrand2008iteratively,  liu2012adaptive, el2010weighted, gelb2019reducing}, which changes the balance between the fidelity and the regularization terms over components so that the method locally controls the regularization effect.}
High-order TV \cite{archibald2016image,  lefkimmiatis2011hessian,chan2000high,  bredies2010total,  setzer2011infimal} takes $n$-th order derivatives into account to obtain piece-wise constant ($n-1$)-th order derivatives. A jump discontinuity
becomes large for a high-order derivative compared to a low-order derivative, and thus the penalty can focus more on discontinuity regions by using a high-order derivative. The order of derivatives can be constant \cite{archibald2016image,  lefkimmiatis2011hessian} or of various orders \cite{chan2000high,  bredies2010total,  setzer2011infimal}.
An inhomogeneous regularization is another class of regularization to recover signals with mixed features. The inhomogeneous regularization uses a non-constant $p$ for the $\ell_p$ regularization to seamlessly utilize the characteristics of various $p$ values,  which was first proposed in \cite{blomgren1997extensions} for image denoising application.  Compared to the weighted TV, the inhomogeneous regularization changes the shape of the regularization constraint, while the weighted TV changes only the scale, not the geometry of the regularization.

There are strategies to assign the spatially varying exponent $p$ over the signal using different references. The method used in \cite{blomgren1997extensions} designs $p$ as a function of the pixel gradient value where it is close to 1 for a large gradient pixel to simulate the TV regularization at sharp edges, and 2 for a small gradient pixel to behave as Tikhonov regularization at smooth regions.  Another method proposed in \cite{chen2010adaptive} assigns $p$ as a function of the second derivative called the difference curvature distinguishing the edge from flat and ramp regions.
In the design of the inhomogeneous exponent distribution $p$, the strategies using the first and the second derivatives of the true signal are sensitive to the quality of the derivative estimation. In the denoising application with a large signal-to-noise ratio (SNR), the measurement data provides a robust estimation for the derivative. However, in a wide range of applications, for example, remote sensing, such high-quality information of the true signal is often unavailable, and thus the design of the inhomogeneous regularization exponent remains as a challenge.

In this work, we propose a method to design a spatially varying $p$ for the inhomogeneous regularization from the incomplete and indirect measurement, such as partial incomplete Fourier measurement data. The proposed method utilizes the statistical information of the true signal. From a single measurement, the method generates a set of samples from a standard reconstruction method. Each reconstruction is not accurate, but the set of samples provides robust statistical information to design the inhomogeneous exponent distribution. A similar idea has been used in the variance-based joint sparsity (VBJS) \cite{adcock2019joint} where point-wise variance statistics is utilized to estimate the weights of the weighted $\ell_1$ regularization. \rev{Also,  the variance information is utilized under the framework of Bayesian inference following the uncertainty reduction of the true signal recovery  \cite{zhang2021empirical}.} Our method is different from other methods in that patches are used to extract statistical information instead of point-wise or pixel-wise information. Also, the proposed method uses interrelation on a patch neighborhood in addition to the point-wise gradient value. This strategy improves the stability and reduces the uncertainty in estimating the local characteristics of the signal that assist in classifying each patch to predefined features with corresponding $p$ values.

In the assignment of $p$ in each patch, the range of $p$ is restricted to $[1,2]$ as in \cite{blomgren1997extensions}. The range of $p$ could extend to $p \in [0,1)$ or $p \in (2,\infty)$ for increased flexibility in the design of the inhomogeneous regularization exponent. $\ell_p$-regularization,  $p \in [0,1)$,  is known to be  superior to $\ell_1$-regularization in capturing sparsity features.  $\ell_p$-regularization,  $\ell_p$, $p \in (2, \infty)$,  is satisfactory to recover very smooth regions.  However,  $\ell_p$,  $p \in [0,1)$ leads to the non-convexity of the objective function and thus requires a delicate optimization method. Also, $\ell_p$, $p \in [2,\infty)$ can suffer from numerical instability while the gain is marginal compared to $p=2$. For this consideration, we restrict $p$ to $[1,2]$, maintaining a convex optimization framework.

The rest of the paper is organized as follows: Section \ref{sec:inhomo} details the problem setup and briefly reviews the inhomogeneous regularization, which is $\ell_p$ with a spatially varying $p$. In section \ref{sec:powerdist}, we propose an intrusive method to design the spatially varying $p$ distribution from indirect and limited data. In Section \ref{sec:admm}, we modify ADMM \cite{boyd2011distributed} as an efficient optimization method to solve the inhomogeneous regularization problem. Section \ref{sec:results} provides 1D and 2D numerical experiments validating the effectiveness and robustness of the proposed method, including recovery of sea ice from incomplete Fourier measurement data. Finally, we conclude this paper with discussions about the limitation and future directions of the current study in Section \ref{sec:conclusions}.

\section{Inhomogeneous regularization}\label{sec:inhomo}
In this section, we formulate the inverse problem of our interest and briefly review the standard regularization method using a penalty term, such as $\ell_p$ and $TV$ regularizations. When the signal has mixed features, such as sparsity and oscillatory behaviors, a homogeneous type regularization can suffer from inaccuracy due to inappropriate characterization of signal characteristics that change over the domain. As a strategy to reconstruct such signals with mixed features, we review the inhomogeneous regularization \cite{blomgren1997extensions}, which is an $\ell_p$ with a spatially varying $\bm{p}$. We also discuss how to design the distribution of $\bm{p}$ for denoising applications and its limitations in other applications.

\subsection{Problem setup}\label{subsec:problem_setup}
We are interested in recovering a signal $\bm{u} \in \mathbb{R}^{N}$ from a measurement $\bm{y} \in \mathbb{R}^{m}$. In particular, we assume that $\bm{u}$ and $\bm{y}$ satisfy the following relation
\begin{equation}\label{eq:linInversProblem}
\bm{y}=\bm{Au}+\bm{\epsilon}.
\end{equation}
Here $\bm{A} \in \mathbb{R}^{m \times N}$ is a linear operator. $\bm{\epsilon}\in \mathbb{R}^{m}$ is a measurement error vector that is assumed to be Gaussian with zero mean. We further assume that the measurement error is uncorrelated across different components, and thus $\bm{\epsilon}$ has a diagonal covariance.
The measurement operator $\bm{A}$ can describe a wide range of measurement operators, such as blurring using a Gaussian kernel. In the case of an image denoising application, $\bm{A}$ is the identity matrix, and the inverse problem is interested in recovering $\bm{u}$ from $\bm{u}+\bm{\epsilon}$. In this study, we are particularly interested in the discrete Fourier measurement that finds applications in magnetic resonance  (MRI) or Synthetic Aperture Radar (SAR) for sea ice imaging.

Under the aforementioned setup, the reconstruction of $\bm{u}$ from $\bm{y}$ is obtained by solving the following optimization problem
\begin{equation}\label{eq:optpb}
\operatorname*{argmin}_{\bm{u} \in \mathbb{R}^N}  \|\bm{A}\bm{u}-\bm{y}\|^2_2.
\end{equation}
When the measurement data is limited or incomplete, there is insufficient information to perfectly reconstruct $\bm{u}$. That is, the reconstruction of $\bm{u}$ from $\bm{y}$ is ill-conditioned when $m<N$. Regularization plays an essential role for the well-posedness of the inverse problem. The regularization forces constraints or restricts $\bm{u}$ to a particular space reflecting prior information. The inverse problem is typically regularized by adding a penalty term $\mathcal{R}(\bm{u})$
\begin{equation}\label{eq:optimization}
\hat{\bm{u}}(\lambda) =  \operatorname*{argmin}_{\bm{u} \in \mathbb{R}^N}  \|\bm{A}\bm{u}-\bm{y}\|^2_2 + \lambda \mathcal{R}(\bm{u}),
\end{equation}
where the regularization parameter $\lambda$ determines the balance between the penalty and the fidelity terms. Tikhonov regularization uses the $\ell_2$ norm of $\bm{u}$ to recover a smoothly varying signal, while $\ell_1$ regularization is useful to reconstruct a signal with sparse structures. The sparsity can be extended to the gradient space of $\bm{u}$; the total variation (TV) regularization, which uses the $\ell_1$ norm of the gradient, can reconstruct a signal with sharp edges or discontinuities.

\subsection{Inhomogeneous regularization using spatially varying $p$}\label{subsec:proposed_regularization}
When the signal has mixed features, for example, edges and oscillations, an application of $\ell_p$ with a homogeneous $p$ regularization will suffer from performance degradation. For example, TV regularization for a smooth oscillatory signal will generate staircase effects while $\ell_2$ regularization fails to capture sparse/discontinuous feature.
To address such issues, \cite{blomgren1997extensions} proposed an inhomogeneous regularization, which is an $\ell_p$ with a spatially varying $l_p$
\begin{equation}\label{proposedRegForm}
\mathcal{R}(\bm{u};\bm{p}): =  \sum\limits_{i=1}^{N}  |Du_i|^{p_i},
\end{equation}
where $D$ is a discrete differential operator, for example, the forward Euler discrete gradient operator (which can be  read as 1-dimensional or 2-dimensional depending on the unknown signal).
Instead of the homogeneous $p_i$ in the standard TV regularization, the inhomogeneous regularization uses the exponent vector $\bm{p}=(p_1,p_2,\cdots, p_{N})$ to adaptively use the characteristics of $\ell_p$ for various $p$ values, while each $p_i, i=1,2,...,N$, is bounded in $[1,2]$ due to the computational benefits discussed in the previous section. \cite{blomgren1997extensions} designed the exponent $p_i$ as a function of the pixel gradient value, $Du_i$, (i.e.,  $p_i=p(Du_i)$) for image denoising application. As the function satisfies $\lim \limits_{x\rightarrow 0}p(x)=2$ and $\lim \limits_{x\rightarrow \infty}p(x)=1$, it adaptively simulates the TV regularization at sharp gradients and Tikhonov regularization at smoother regions.

We note that the inhomogeneous regularization is different from the weighted TV
\begin{equation}\label{eq:wTV}
\mathcal{R}(\bm{u})=\sum_{i=1}^{N}w_i|Du_i|
\end{equation}
using a weight vector $\bm{w}\geq \bm{0}$. The weighted TV changes the scale of each $|Du_i|$, while preserving the geometry of the constraint. On the other hand, the inhomogeneous regularization changes the geometry of the constraint while maintaining the convexity. It is natural to combine the two methods to utilize the characteristics of each method
 \begin{equation}\label{curvatureRegForm}
  \sum\limits_{i=1}^{N} w_i |Du_i|^{p_i}
\end{equation}
using an exponent vector $\bm{p}$ and a weight vector $\bm{w}$. A PDE-based denoising method \cite{chen2010adaptive} considered both the spatially varying $\ell_p$ and the balance, in which pointwise weight is applied to the corresponding fidelity term compared to Eq.~\eqref{curvatureRegForm}.  Exponents and weights are designed as functions of the second derivative based edge indicator called difference curvature.  The edge indicator is the difference between second directional derivatives in tangential and normal directions of the gradient,  which has better performance than the gradient and other curvatures in distinguishing edges from flat and ramp regions and isolated noise.  The authors assign exponents and weights proportional to the square root of the difference curvature following the similar idea in \cite{blomgren1997extensions} mentioned above.

The focus of the current study is the application of inhomogeneous regularization to the case of incomplete and indirect data. The aforementioned strategies to design the exponent vector $\bm{p}$ requires the first and/or second order derivatives of the true signal. In the denoising problem, as the measurement operator is the identity, such information of the unknown signal can be estimated when the signal-to-noise is low. However, if the measurement operator is limited and indirect, no such information is available. In particular, there is no robust exponent design strategy for an incomplete measurement of the Fourier data.

\section{Construction of the exponent distribution in the inhomogeneous regularization} \label{sec:powerdist}
In this section, we propose a strategy to design the exponent distribution $\bm{p}$ of the inhomogeneous regularization \eqref{proposedRegForm} under incomplete indirect data. Our numerical experiments in the next section show the robust performance of the inhomogeneous regularization without weights, and thus we focus on the inhomogeneous exponent in the current study, while we leave the weight design as future work.
In the design of the inhomogeneous exponent following the ideas in \cite{blomgren1997extensions,chen2010adaptive}, the main issue of the indirect and limited measurement data  is the lack of an accurate estimate of the first and/or the second derivative information. We address this issue by estimating the exponent based on statistical information of the unknown signal. In particular, we generate samples using multiple reconstructions from the standard regularization, such as TV or $\ell_2$. Any single reconstruction using them is not necessarily accurate as the regularization can be inappropriate. However, using patches instead of pointwise values, we can extract robust statistical information to design the inhomogeneous exponent.
That is, we identify the local characteristic of the signal within a small patch rather than a single point.
The rationale of the patches is that the signal characteristics to design $\bm{p}$, such as sparsity and oscillation, rely on the interrelation with neighborhoods.

Let $\{\Omega_j\}_{j=1}^{M}$ be the set of non-overlapping patches that covers the domain of the signal. Then, the regularization functional \eqref{proposedRegForm} can be rewritten as
\begin{equation}\label{patchwiseProposedRegForm}
\mathcal{R}(\bm{u};\bm{p}) =  \sum\limits_{j=1}^{M}\sum\limits_{i \in \mathcal{I}_j}  |Du_i|^{q_j},
\end{equation}
where $\mathcal{I}_j$ is the index set corresponding to the patch $\Omega_j$, and $q_j$ is an exponent that is constant in $\Omega_j$.
Here, we presume that the regularization effect is local within a patch as inter-patch interactions do not diffuse very far \cite{louchet2011total}.
For a 2D image, we can consider the total variation $|Du_i|$ as either isotropic TV, $|Du_i|=\sqrt{(\partial_x u_i)^2+(\partial_y u_i)^2}$, or the anisotropic TV, $|Du_i|=|\partial_xu_i|+|\partial_yu_i|$. One difference between the two candidates is the rotational invariant property of the isotropic TV, and each type is employed in favor of applications. In this work, we consider the isotropic TV for its convenience in the adaptation in ADMM, which is explained in Section \ref{sec:admm}. Our construction of the exponent distribution consists of three sub-steps; i) generating multiple samples from a single measurement data, ii) classifying each patch to predefined features, and iii) assigning the exponent value based on the classification.

\subsection{Multiple reconstructions from a single data}
From a single measurement data $\bm{y}$, we reconstruct multiple samples by changing the balance between the fidelity and the homogeneous regularization penalty terms. That is, we solve the following standard homogeneous regularization using different values of $\lambda$
\begin{equation}\label{eq:optProposedReg}
\hat{\bm{u}}(\lambda;\bm{p}) :=  \operatorname*{argmin}_{\bm{u} \in \mathbb{R}^N} \left(  \|\bm{A}\bm{u}-\bm{y}\|_2^2 + \lambda \mathcal{R}(\bm{u};\bm{p}) \right),
\end{equation}
where $\bm{p}$ is a constant vector. This idea of generating multiple measurement samples has been used in \cite{adcock2019joint} to design the weight of the weighted $\ell_1$ regularization, where only one value of $\bm{p}$ has been used. In our study, we use variance-based statistics to extract features selectively from two sets of samples; discontinuity from TV-regularized samples ($\bm{p}=\bm{1}$) and oscillation from Tikhonov-regularized samples ($\bm{p}=\bm{2}$). For each case, we use a logarithmically-spaced distribution to draw $\lambda$. The rationale behind the logarithmically-spaced distribution is to have samples not necessarily biased by either the fidelity or the regularization terms. By drawing multiple $\lambda$ values, we solve the standard optimization \eqref{eq:optProposedReg}, which yields two sets of reconstructed samples, $\{\hat{\bm{u}}(\lambda_i;\bm{1})\}_{i=1}^{C}$ and $\{\hat{\bm{u}}(\lambda_i;\bm{2})\}_{i=2}^{C}$, from the TV and Tikhonov regularizations, respectively.

\subsection{Patch classification}
From the multiple reconstructions, we classify each patch into three categories, i) discontinuity, ii) oscillation, and iii) smoothness.
First, we estimate the gradient of the unknown signal using the averaged gradient of the reconstructed samples $\{\hat{\bm{u}}(\lambda_i;\bm{1})\}_{i=1}^{C}$ and ${\hat{\bm{u}}(\lambda_i;\bm{2})\}_{i=2}^{C}}$
\begin{equation}
\bm{g}_1=(g_{1,1},g_{1,2},\cdots, g_{1,N}) = \frac{1}{C} \sum \limits_{i=1}^{C}D\hat{\bm{u}}(\lambda_i;\bm{1}),
\end{equation}
\begin{equation}
\bm{g}_2=(g_{2,1},g_{2,2},\cdots, g_{2,N}) = \frac{1}{C} \sum \limits_{i=1}^{C}D\hat{\bm{u}}(\lambda_i;\bm{2}).
\end{equation}
To classify the patches in an unsupervised way, we compute the patch-wise variance of the gradient called variance pooling
\begin{equation}
\operatorname*{var}(\bm{g}_l;\Omega_j) = \frac{1}{|\mathcal{I}_j|}\sum \limits_{i \in \mathcal{I}_j}g_{l,i}^2-\left(\frac{1}{|\mathcal{I}_j|}\sum \limits_{i \in \mathcal{I}_j}\left |g_{l,i}\right| \right)^2, ~l=1,2,
\end{equation}
\begin{algorithm}[H]
 \caption{Local feature classification \label{algo:patchClassification}}
\begin{algorithmic}[1]
\Require{patch size $K$,  variance threshold $\epsilon$,  MinMax-neighborhood filter size $m$}
\Statex
\State Choose the set of regularization parameters $\lambda_i=1,\cdots, C$.  Recover the unknown signal from TV and Tikhonov regularizations;

$$\hat{\bm{u}}(\lambda_i;\bm{p}) =  \operatorname*{argmin}_{\bm{u} \in \mathbb{R}^N} \left(  \|\bm{A}\bm{u}-\bm{y}\|_2^2 + \lambda_i \mathcal{R}(\bm{u};\bm{p}) \right),~~\bm{p}=\bm{1},\bm{2}.$$
\\
Compute the gradient statistics from $\{\hat{\bm{u}}(\lambda_i;\bm{1})\}_{i=1}^{C}$ and $\{\hat{\bm{u}}(\lambda_i;\bm{2})\}_{i=1}^{C}$;

$$\bm{g}_1 = \frac{1}{C} \sum \limits_{i=1}^{C}D\hat{\bm{u}}(\lambda_i;\bm{1}), ~~~\bm{g}_2 = \frac{1}{C} \sum \limits_{i=1}^{C}D\hat{\bm{u}}(\lambda_i;\bm{2}).$$
\\
Split the gradient map $\bm{g}_1$ and $\bm{g}_2$  to $K$-sized patches,  $\{\Omega_j\}_{j=1}^{M}$.  Apply the variance pooling in each patch,  $\Omega_j$;

$$\operatorname*{var}(\bm{g}_l;\Omega_j) = \frac{1}{|\mathcal{I}_j|}\sum \limits_{i \in \mathcal{I}_j}g_{l,i}^2-\left(\frac{1}{|\mathcal{I}_j|}\sum \limits_{i \in \mathcal{I}_j}\left |g_{l,i}\right| \right)^2,~~l=1,2,$$
where $\mathcal{I}_j$ is the index set of the patch $\Omega_j$.
\\
Classify the patches $\{\Omega_j\}_{j=1}^{M}$ into $3$ categories;
$$\operatorname*{class}(\Omega_j)=
\begin{cases}
\text{discontinuity} & \textrm{~if~} (\bm{g}_1\text{-condition}) \text{ and } \textbf{not}(\bm{g}_2\text{-condition}), \\
\text{oscillation} & \textrm{~if~} (\bm{g}_2\text{-condition}) \text{ and } \textbf{not}(\bm{g}_1\text{-condition}), \\
\text{smooth} & \textrm{~otherwise} , \\
\end{cases}$$
where
$$
\begin{cases}
\bm{g}_1\text{-condition}= \widetilde{\operatorname*{var}}(\bm{g}_1;\Omega_j) \geq \epsilon \textrm{~and~} \operatorname*{NghdFilter}(\bm{g}_1, \Omega_j,  m) < \epsilon, \\
\bm{g}_2\text{-condition}=\widetilde{\operatorname*{var}}(\bm{g}_2;\Omega_j) \geq \epsilon \textrm{~and~} \operatorname*{NghdFilter}(\bm{g}_2, \Omega_j,  m) \geq \epsilon,
\end{cases}
$$
and $\widetilde{\operatorname*{var}}(\bm{g}_l;\Omega_j)$ is the Min-Max normalization of $\operatorname*{var}(\bm{g}_l;\Omega_j)$, and $\operatorname*{NghdFilter}(\Omega_j,  m)$ is read in Section \ref{sec:powerdist}.
\end{algorithmic}
\end{algorithm}
\noindent
and normalize the variance map with Min-Max normalization (denoted as $\widetilde{\operatorname*{var}}(\bm{g}_l;\Omega_j)$).  The variance is expected to be low if a patch lies on a smooth region and high if there is discontinuity or oscillation. We identify the smooth region by introducing the threshold value $\epsilon$ on the variance. The patch with a variance greater than the threshold $\epsilon$ is classified as either discontinuity or oscillation.
We distinguish the discontinuity patch from the oscillation patch by observing the neighborhood patches.

We first define the $n$-sized directional neighborhoods for given patch $\Omega_{i}$ (1D) or $\Omega_{ij}$ (2D) as
\begin{eqnarray}
\mbox{(1D)}\quad \mathcal{N}&=&\bigcup\limits_{k=1}^{n} \Omega_{i\pm k},  \\
\mbox{(2D)}\quad \mathcal{N}_1&=&\bigcup\limits_{k=1}^{n} \Omega_{i,j\pm k}, ~\mathcal{N}_2=\bigcup\limits_{k=1}^{n} \Omega_{i\pm k,j}, ~ \mathcal{N}_3=\bigcup\limits_{k=1}^{n} \Omega_{i\pm k,j \pm k}, ~\mathcal{N}_4=\bigcup\limits_{k=1}^{n} \Omega_{i\mp k,j\pm k}
\end{eqnarray}
Note that the centered patch is exclusive in each neighborhood, and four directional neighborhoods (i.e., horizontal, vertical, and two diagonal directions) are considered in a 2D image. We utilize the different characteristics of two classes in the neighborhoods; a discontinuity patch is isolated from the surrounding smooth patches. On the other hand, an oscillatory patch is encompassed by other oscillatory patches. Thus, a discontinuity patch has a directional patch-neighborhood with low variance patches, which is not the case for oscillation. To capture such behavioral difference, we develop a filter on the variance-pooled gradient map,
\begin{equation}
\operatorname*{NghdFilter}(\bm{g}_l, \Omega_i, n) = \max \limits_{1\leq k \leq n} v^{(l)}_{i \pm k} \tag{\text{1D}}
\end{equation}
\begin{equation}
\operatorname*{NghdFilter}(\bm{g}_l,\Omega_{i,j}, n)=
\min \left\{\max \limits_{1\leq k \leq n} v^{(l)}_{i, j \pm k},
\max \limits_{1\leq k \leq n} v^{(l)}_{i \pm k, j},
\max \limits_{1\leq k \leq n} v^{(l)}_{i,\pm k j \pm k},
\max \limits_{1\leq k \leq n} v^{(l)}_{i,\mp k j \pm k}
\right\} \tag{\text{2D}}
\end{equation}
which computes the minimum among the maximums of each directional neighborhood. Here $ v^{(l)}_{p} = \widetilde{ \operatorname*{var}}(\bm{g}_l;\Omega_{p})$ and  $ v^{(l)}_{p,q} = \widetilde{ \operatorname*{var}}(\bm{g}_l;\Omega_{p,q})$.

A filtered value less than the smooth threshold $\epsilon$ used above implies that there is a smooth region around, and thus the patch is classified as discontinuity.
Otherwise, non-smooth regions surround the patch, and thus the patch is identified as oscillation.
To identify discontinuity, $\bm{g}_1$-statistics are preferred to $\bm{g}_2$-statistics
 in that the first have the higher value on the discontinuous patch and the smooth neighborhoods are well-identified with lower value since the second is contaminated by the Gibb's phenomena around the discontinuity.  On the other hand,  $\bm{g}_2$-statistics are more stable for identifying oscillation since the $\bm{g}_1$-statistics have a chance to mislead the oscillatory feature by the staircase effect (a numerical test for  comparison between $\bm{g}_1$ and $\bm{g}_2$ statistics in variance pooling maps are presented in Fig.~\ref{fig:synthetic_1D_pooling}(a) and (b) in Section~\ref{sec:results}). The complete local feature classification algorithm is summarized in Algorithm~\ref{algo:patchClassification}.

\subsection{Exponent distribution}
We interpret $\ell_p, 1<p<2$, as an interpolation balancing between $\ell_1$ and $\ell_2$ regularizations. 
The sparsity is relaxed as the exponent is far away from $1$, and the regularization weights more on smoothness as the exponent gets close to $2$ in the sense of Tikhonov regularization.
Following this interpretation, we assign $1$ for the exponent on a discontinuity patch.  For smooth and oscillatory patches,  we adaptively set the exponent in reference to the patch-wise average of the gradient (called average pooling),

\begin{algorithm}[H]
 \caption{Signal recovery with the inhomogenous regularization \label{algo:wholeProcess}}
\begin{algorithmic}[1]
\Require{patch size $K$,  variance threshold $\epsilon$,  MinMax-neighborhood filter size $m$, power distribution constant $c$}
\Statex
\State Choose the set of regularization parameters $\lambda_i=1,\cdots, C$.  Recover the unknown signal from TV and Tikhonov regularizations;

$$\hat{\bm{u}}(\lambda_i;\bm{p}) =  \operatorname*{argmin}_{\bm{u} \in \mathbb{R}^N} \left(  \|\bm{A}\bm{u}-\bm{y}\|_2^2 + \lambda_i \mathcal{R}(\bm{u};\bm{p}) \right),~~\bm{p}=\bm{1},\bm{2}$$
\\
Compute the gradient statistics from $\{\hat{\bm{u}}(\lambda_i;\bm{1})\}_{i=1}^{C}$ and $\{\hat{\bm{u}}(\lambda_i;\bm{2})\}_{i=1}^{C}$;

$$\bm{g}_1 = \frac{1}{C} \sum \limits_{i=1}^{C}D\hat{\bm{u}}(\lambda_i;\bm{1}), ~~~\bm{g}_2 = \frac{1}{C} \sum \limits_{i=1}^{C}D\hat{\bm{u}}(\lambda_i;\bm{2})$$
\\
Split the gradient map $\bm{g}_1$ and $\bm{g}_2$  to $K$-sized patches,  $\{\Omega_j\}_{j=1}^{M}$.   Classify patches with Algorithm~\ref{algo:patchClassification}.
\\ Apply the average pooling on $\bm{g}_1$ and $\bm{g}_2$ maps;
$$\operatorname*{avg}(\bm{g}_l;\Omega_j) = \frac{1}{|\mathcal{I}_j|}\sum \limits_{i \in \mathcal{I}_j}|g_{l,i}|,~l=1,2$$
where $\mathcal{I}_j$ is the index set of the patch $\Omega_j$.
\\ Distribute the exponent on each patch;
$$
\Phi(\Omega_j) =
\begin{cases}
1 & \textrm{~if~} \Omega_j \in \textrm{discontinuity class}, \\
2 - \exp \left(-c\cdot \widetilde{\operatorname*{avg}}(\bm{g}_1;\Omega_j) \right)& \textrm{~if~} \Omega_j \in  \textrm{smooth class},\\
2 - \exp \left(-c\cdot \widetilde{\operatorname*{avg}}(\bm{g}_2;\Omega_j) \right)& \textrm{~if~} \Omega_j \in \textrm{oscillation class}, \\
\end{cases}
$$
where $ \widetilde{\operatorname*{avg}}(\bm{g}_l;\Omega_j) $ is Min-Max normalization of $\operatorname*{avg}(\bm{g}_l;\Omega_j)$.
\\  Recover the unknown signal with inhomogeneous exponent distribution $\{\Phi(\Omega_j)\}_{j=1}^{M}$,
$$
\hat{\bm{u}}\left(\lambda_l;\{\Phi(\Omega_j)\}_{j=1}^{M}\right)=\operatorname*{argmin}_{\bm{u} \in \mathbb{R}^N} \left(  \|\bm{A}\bm{u}-\bm{y}\|_2^2 +  \lambda_l\sum\limits_{j=1}^{M}\sum\limits_{i \in \mathcal{I}_j}  |Du_i|^{\Phi(\Omega_j)}\right).
$$
\end{algorithmic}
\end{algorithm}

\begin{equation}
\operatorname*{avg}(\bm{g}_l;\Omega_j) = \frac{1}{|\mathcal{I}_j|}\sum \limits_{i \in \mathcal{I}_j}|g_{l,i}|, ~~l=1,2.
\end{equation}
The distribution of the exponent is given as a function $\Phi$ of Min-Max normalization of the average pooling, (denoted as $\widetilde{\operatorname*{avg}}(\bm{g}_l;\Omega_j)$),  which is designed to be an increasing function as a larger exponent induces less penalization on the gradient. $\Phi$ sets the exponent close to  $1$ 
on the flat region (i.e., zero gradient), while it sets the exponent close to $2$ on the oscillatory patch with a high averaged gradient. We design $\Phi$ to be
\begin{equation}
\Phi(\widetilde{\operatorname*{avg}}(\bm{g}_l;\Omega_j)) = 2 - \exp \left(-c\cdot \widetilde{\operatorname*{avg}}(\bm{g}_l;\Omega_j) \right),  ~~ c>0.
\end{equation}

In the numerical experiment in Section~\ref{sec:results} (Fig.~\ref{fig:synthetic_1D_pooling} (c) and (d)),  we observe that on smooth patches,  the $\bm{g}_1$-average pooling is more stable under the choice of reconstructions in comparison with the $\bm{g}_2$-average pooling (and vice versa on oscillatory patches).
Accordingly, we assign $\Phi(\widetilde{\operatorname*{avg}}(\bm{g}_1;\Omega_j))$ on a smooth patch and $\Phi(\widetilde{\operatorname*{avg}}(\bm{g}_2;\Omega_j))$ on an oscillatory patch. The complete procedure is presented in Algorithm~\ref{algo:wholeProcess}.

\section{ADMM for inhomogeneous regularization}\label{sec:admm}
As the exponents of the inhomogeneous regularization remain bounded in $[1,2]$,  a convex optimization method can be utilized to solve the problem.  Among other convex optimization methods, we find that the inhomogeneous regularization can be efficiently solved using the alternating direction method of multipliers (ADMM; \cite{boyd2011distributed}) with an appropriate problem transformation.  ADMM is a variant of the augmented Lagrangian method, which is flexible in decomposing the original optimization into small local subproblems with computational benefits.  For instance,  it reduces the complexity by avoiding the joint (sub)optimizations, or is suitable for parallel computing.  We take such advantage of ADMM to solve the regularized optimization by decoupling the regularization term from the fidelity term and handling both separately.  Among many possibilities to transform the problem to fit in the framework,  we design the transformation to make the following subproblems efficiently solved.  Motivated by the group Lasso regularization \cite{yuan2006model},  we reduce the suboptimization corresponding to the inhomogeneous regularization separable,  which is computationally beneficial.

We rewrite the problem \eqref{eq:optProposedReg} by introducing the auxiliary variable $\bm{v}$ and the constraint as following;
\begin{align}
\text{minimize} &~~ \mathcal{Q}(\bm{u}) + \lambda \mathcal{R}(\bm{v}), \label{eq:ADMM Objectives}\\
\text{subject to} &~~ \bm{F}\bm{u} - \bm{v} = \bm{0}, \label{eq:ADMM constraint}
\end{align}
where $\bm{F} : \mathbb{R}^{N} \rightarrow \mathbb{R}^{M}$ is the discrete gradient operator, $M$ is either $N$ or $2N$ depending on whether the unknown $\bm{u}$ is 1D signal or 2D image,  $\bm{v} =(\bm{v}_1, \bm{v}_2,\cdots, \bm{v}_N)\in \mathbb{R}^{M}$ is the auxiliary variable,  and the objective functions $\mathcal{Q}$ and $\mathcal{R}$ are
\begin{align}
\mathcal{Q}(\bm{u}) &=  \| \bm{A}\bm{u}- \bm{y} \|^2_2, \\
\mathcal{R}(\bm{v}) &= \sum \limits_{i=1}^{N}g_i(\|\bm{v}_i\|_2), ~~ g_i(z) = |z|^{p_i}.\label{eq:rerrange_reg}
\end{align}
We note that the constraint Eq.~\eqref{eq:ADMM constraint} is  read as $\bm{v}_i =Du_i \in \mathbb{R},  i=1,2,\cdots, N$,  for the case of 1D signal,  and $\bm{v}_i = (D_xu_i, D_yu_i) \in \mathbb{R}^2,  i=1,2,\cdots, N$,  for the case of 2D image.
The augmented Lagrangian corresponding to the optimization Eq.~\eqref{eq:ADMM Objectives} and Eq.~\eqref{eq:ADMM constraint} is
\begin{equation}
L_{\rho}(\bm{u}, \bm{v}, \bm{w}) = \mathcal{Q}(\bm{u}) + \rev{\lambda}\mathcal{R}(\bm{v}) + \bm{w}^{T}(\bm{F}\bm{u} - \bm{v}) + (\rho/2)\|\bm{F}\bm{u} - \bm{v}  \|_2^2,
\end{equation}
where $\bm{w}\in \mathbb{R}^{M}$ is the Lagrange multiplier and $\rho>0$ is the penalty parameter.
ADMM consists of the iterations with dual updates
\begin{equation}
\bm{u}^{k+1} =\text{Prox}_{\mathcal{Q},\rho,\bm{F}}(\bm{v}^{k} - \bm{w}^{k}):=\operatorname*{argmin}_{\bm{u}} \left(\mathcal{Q}(\bm{u})+ \frac{\rho}{2}\|\bm{F}\bm{u} -(\bm{v}^{k} - \bm{w}^{k}) \|_2^2\right), \label{eq:u-proximal}
\end{equation}
\begin{equation}
\bm{v}^{k+1} =\text{Prox}_{\mathcal{R},\frac{\rho}{\lambda},\bm{I}}(\bm{F}\bm{u}^{k+1} + \bm{w}^{k}):=\operatorname*{argmin}_{\bm{v}} \left( \mathcal{R}(\bm{v}) + \frac{\rho/\lambda}{2}\| \bm{I}\bm{v} - (\bm{F}\bm{u}^{k+1}  + \bm{w}^{k}) \|_2^2\right),\label{eq:v-proximal}
\end{equation}
\begin{equation}
\bm{w}^{k+1} = \bm{w}^{k}+ \bm{F}\bm{u}^{k+1} - \bm{v}^{k+1}.
\end{equation}
The suboptimizations,  Eq.~\eqref{eq:u-proximal} and Eq.~\eqref{eq:v-proximal} are known as proximal operators, which,  in general,  have either closed-form or implicit form depending on the objective terms.  The $\bm{u}$-optimization Eq.~\eqref{eq:u-proximal},  as a quadratic objective form,  is known to have the  closed-form
\begin{equation}
\text{Prox}_{\mathcal{Q},\rho,\bm{F}}(\bm{x})= (\bm{A}^T\bm{A}+\rho\bm{F}^{T}\bm{F})^{-1}(\bm{A}^T\bm{b}+\rho\bm{F}^{T}\bm{x}), ~~\bm{x} \in \mathbb{R}^{M}.
\end{equation}
For the proximal operator in $\bm{v}$-optimization Eq.~\eqref{eq:v-proximal},  the objective function $\mathcal{R}(\bm{v})$ is separable,  so is the optimization.  For $\bm{x}=(\bm{x}_1, \bm{x}_1,\cdots, \bm{x}_N) \in \mathbb{R}^{M}$,
\begin{align}
\left(\text{Prox}_{\mathcal{R},\rho,\bm{I}}(\bm{x})\right)_i&= \operatorname*{argmin}_{\bm{v}_i} \left( g_i(\|\bm{v}_i\|_2) + \frac{\rho}{2}\| \bm{v_i} - \bm{x_i} \|_2^2 \right) \\ \label{eq:prox_R}
&=\text{Prox}_{g_i(\|\cdot \|_2), \rho, \bm{I}}(\bm{x}_i) \\
&= \text{Prox}_{g_i, \rho, \bm{I}}(\|\bm{x}_i\|_2)\frac{\bm{x}_i}{\|\bm{x}_i\|_2},  ~~ i=1,2,\cdots,N, \label{eq:prox_norm_composition}
\end{align}
in which the last equality holds by the norm composition rule of the proximal operator.  We can easily derive the proximal operator,  $\text{Prox}_{g_i, \rho, \bm{I}}(q)$ in Eq.~\eqref{eq:prox_norm_composition} as
\begin{equation}\label{eq:v-update}
\text{Prox}_{g_i, \rho, \bm{I}}(q) =
\begin{cases}
\mathcal{S}_{1/\rho}(q) & ~\text{if } p_i =1,
\\
\text{The zero of }~ h(x) :=  \text{sgn}(x)p_i|x|^{p_i-1} + \rho (x - q) & ~\text{if } p_i > 1,
\end{cases}
\end{equation}
where $\mathcal{S}_{\kappa}$ is the soft thresholding or shrinkage operator defined as $\mathcal{S}_{\kappa}(q):=\left(1-\frac{\kappa}{q} \right)_{+}q$.  We note that $h(x)$ has the unique zero as it is an increasing function.  We find the root of $h(x)$ with the Chadrupatla's method \cite{chandrupatla1997new},  which is superior to find the root of a function with stiff regions,  as $h(x)$ does near $0$,  than secant methods.

The algorithm is also applied for the case of anisotropic TV by replacing $g_i(\|\bm{v}_i\|_2)$ in Eq.~\eqref{eq:rerrange_reg} to $g_i(\|\bm{v}_i\|_1)$.  The $\mathcal{Q}$-proximal operator is invariant,  but the $\mathcal{R}$-proximal operator corresponding to Eq.~\eqref{eq:prox_R} should be directly solved by an iterative method due to the limitation of composition rule for $\ell_1$-norm.

\section{Numerical Results}\label{sec:results}
We validate the robustness and effectiveness of the proposed method to design the inhomogeneous regularization through a suite of numerical tests, including sea ice reconstruction. In particular, as an incomplete and indirect measurement, we consider partial Fourier measurement used in many applications,  such as magnetic resonance (MRI) or Synthetic Aperture Radar (SAR) imaging.  The measurement is characterized as a global operator in that each Fourier coefficient depends on the global shape of the signal,  which contrasts to the local measurement such as blurred or noisy signal.
We conduct a numerical experiment to investigate the consistency of the exponent design under various sets of reconstructed samples. We particularly analyze the distributions of variance and average pooling maps from multiple sets of reconstruction samples, each set of which corresponds to the logarithmically-spaced regularization parameters on a randomly generated interval.
In generating multiple samples, we can easily parallelize to solve multiple optimization problems with different regularization parameters which are not coupled with each other.  In the ADMM framework,  a regularization parameter is only engaged in $\bm{v}$-optimization, Eq.~\eqref{eq:v-proximal}, in the dual update scheme, which is also separable. We compute the $\bm{v}$-update in parallel in both regularization parameters and components by Eq.~\eqref{eq:v-update}.
In addition, we use the penalty parameter $\rho=1$ in ADMM for all examples.

We measure the performance of the proposed method using the point-wise and $\mathcal{L}_p$-errors,  for $p=1,2$,  \rev{comparing with other regularizations; TV,  Tikhonov,  the curvature-based inhomogeneous regularization \cite{chen2010adaptive} where the curvature is calculated using the average of the samples}
\begin{equation}
\frac{1}{C} \sum \limits_{i=1}^{C}\hat{\bm{u}}(\lambda_i;\bm{1}).
\end{equation}
\rev{We also compare the VBJS regularization \cite{adcock2019joint} in the examples of synthetic signal recovery.}
\begin{figure}[h]
\centering
\includegraphics[width=1.0\textwidth]{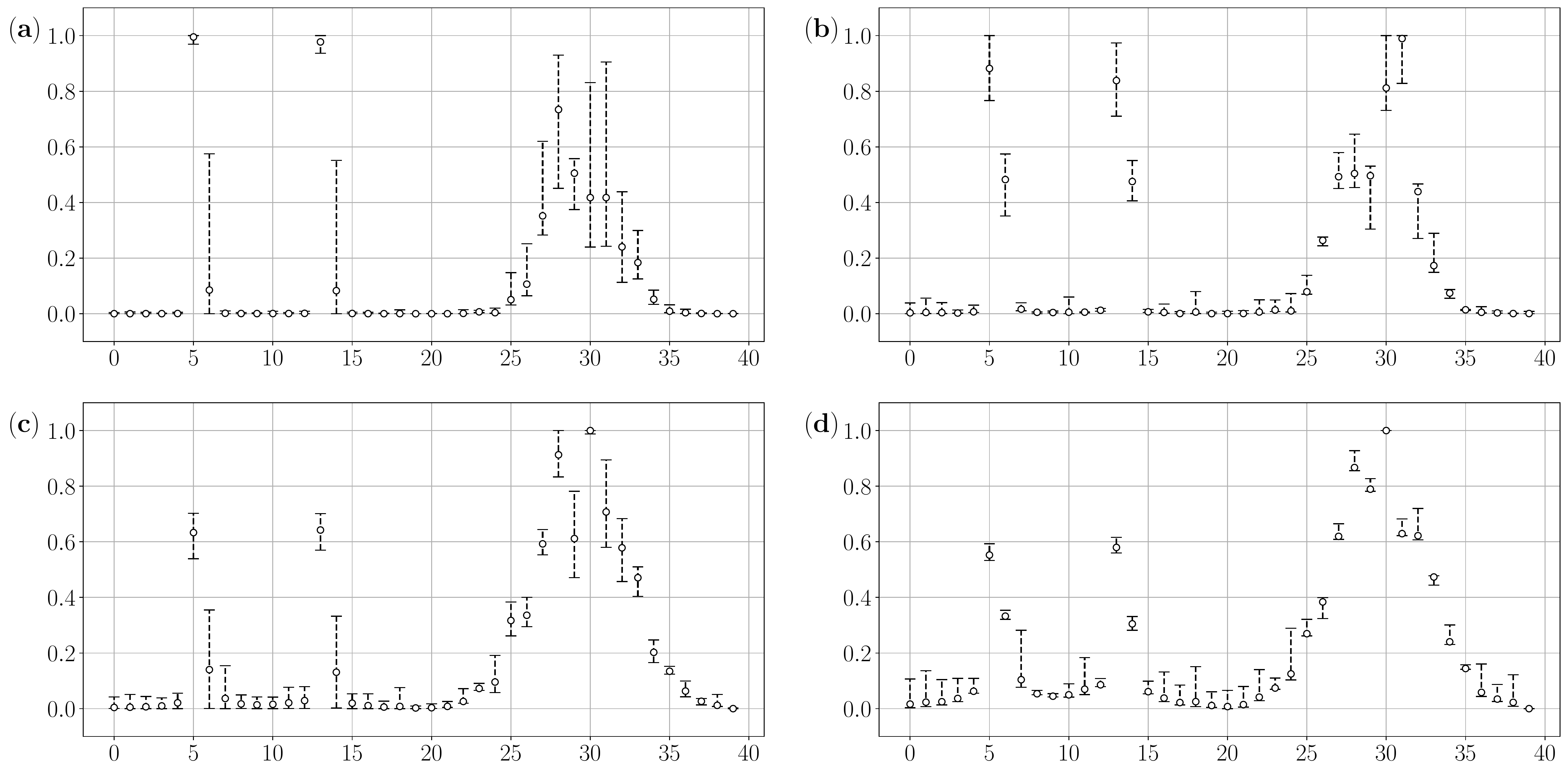}
\caption{The distributions of the pooling maps ($K=5$) from multiple sets of reconstruction samples.  The first column ((a) and (c)) describes the pooling maps from $\{\hat{\bm{u}}(\lambda_i^{(j)},\bm{1})\}_{i=1}^{C}$ and the second column ((b) and (d)) corresponds to $\{\hat{\bm{u}}(\lambda_i^{(j)},\bm{2})\}_{i=1}^{C}$.  The first row ((a) and (b)) and the second row ((c) and (d)) describe the variance and average pooling maps, respectively.  The empty dot represents the mean of pooled samples and the bar corresponds to the range of pooled values on each patch.}
\label{fig:synthetic_1D_pooling}
\end{figure}

\subsection{Synthetic 1D signal recovery}\label{subsec:syn_1d}
The first test is a synthetic 1D signal. The signal is generated by combining two signals with
different characteristics
\begin{equation}
u:[-1,1] \rightarrow \mathbb{R}, \hspace{1mm} u(x) =
\begin{cases}
1 & ~\textrm{if}~ -0.7 \leq x \leq -0.3 \\
\frac{1}{2}(1+\sin(100(x+1))\times \exp(-25(x-0.5)^2) & ~\textrm{if}~ 0 < x \leq 1 \\
0 & ~\textrm{otherwise}
\end{cases}
\end{equation}
where the signal has a sparse gradient feature on $[-1,0]$ and highly oscillatory behavior on $(0,1]$. We discretize the signal with $N=200$ uniform grid points in $[-1,1]$, which yields a 200-dimensional true signal.

The measurement is the first $20\%$ of the smallest wavenumber components of the Fourier transform of the true signal.
As each coefficient is a complex number,  we transform it to the real-valued vectors by separating the real and imaginary parts (i.e., Fourier cosine and sine coefficients); thus, the total number of measurement is $40\%$ of the signal length (i.e., $m=\frac{40}{100}N$). To focus on the design of the inhomogeneous regularization while minimizing the effect from the measurement error, we do not contaminate the measurement with the Gaussian noise in this test.

Fig.~\ref{fig:synthetic_1D_pooling} shows the distribution of variance and average pooling maps generated from 1000 TV and Tikhonov samples, $\{\hat{\bm{u}}(\lambda_i^{(j)},\bm{1})\}_{i=1}^{200}$ and $\{\hat{\bm{u}}(\lambda_i^{(j)},\bm{2})\}_{i=1}^{200}$, $j=1,2,\cdots,1000$.
Each regularization parameter set, $\{\lambda_i^{(j)}\}_{i=1}^{200}$, is logarithmically-spaced values in a random interval from $[10^{-4}, 10^{4}]$ with random length in the range from $10^2$ to $10^4$.
We compare two distributions of variance pooling maps (Fig.~\ref{fig:synthetic_1D_pooling}(a) and (b)) in the aspect of predefined features for the classification.
Smooth regions consistently have low variance in both distributions,  which are distinguished from the discontinuity and oscillatory regions with high variance.  The discontinuity feature is more discernible in TV-regularized reconstructions (Fig.~\ref{fig:synthetic_1D_pooling}(a)) in that a patch with a discontinuity is well-isolated from low variance patches identified as the smooth class,  and is stable in high variance value.  On the other hand,  variance values on oscillatory regions are  higher in average and more consistent  in Tikhonov-regularized 
\begin{figure}[!]
\centering
\includegraphics[width=1\textwidth]{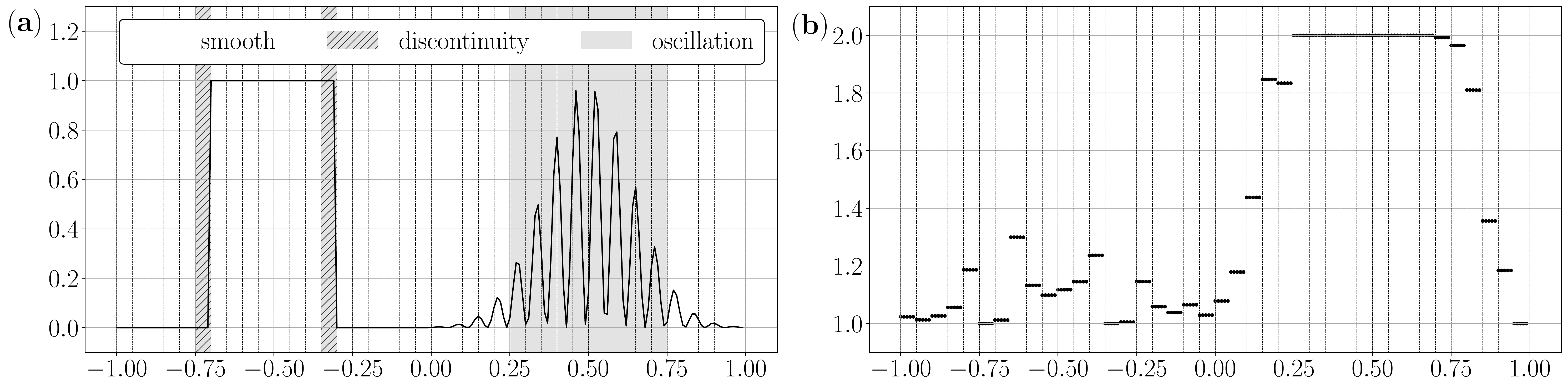}
\caption{Local feature classification (a) and exponent distribution (b) with the hyperparameters $(K,\epsilon, m, c)=(5,1\text{E-}2,3,27)$ in algorithm~\ref{algo:wholeProcess} }
\label{fig:example_1_classification_exponent}
\end{figure}
\begin{figure}[h]
\centering
\includegraphics[width=1.0\textwidth]{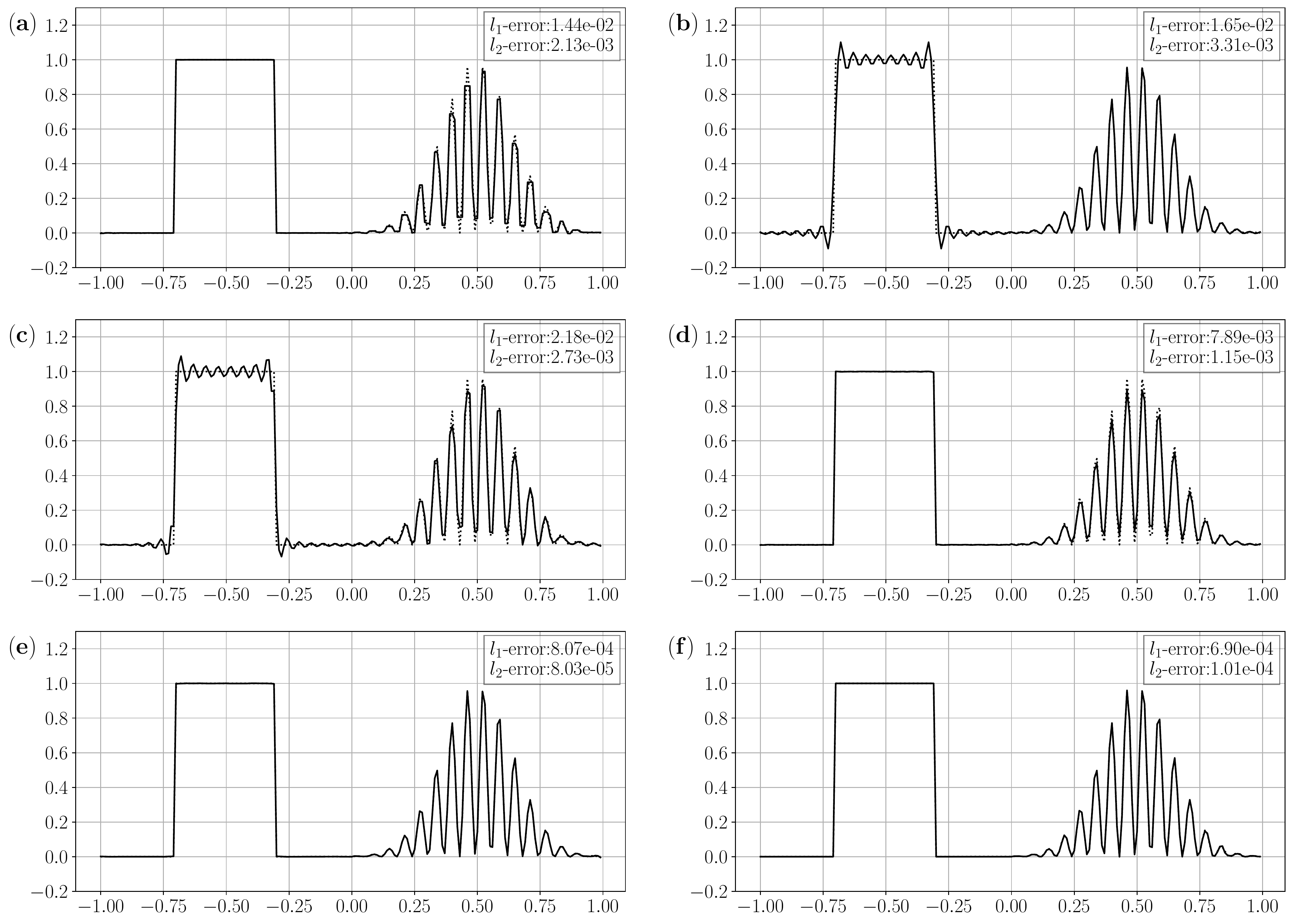}
\caption{\rev{1-dimensional signal recovery from six different regularizations; (a) homogeneous $\bm{p}=\bm{1}$,  (b) homogeneous $\bm{p}=\bm{2}$,  (c) curvature-based regularization,  (d) VBJS regularization,  (e)  $\bm{p}=\bm{1} $ on $[-1,0]$ and $\bm{p}=\bm{2}$ on $(0,1]$,  (f) proposed regularization}}
\label{fig:example_1_approximation}
\end{figure}
reconstructions (Fig.~\ref{fig:synthetic_1D_pooling}(b)).  We note that the identification for smooth regions is a baseline for our feature classification algorithm,
and the variance statistic has more consistent distribution on smooth regions than the average statistic.  Moreover, a  discontinuity  is more distinguishable in both magnitude and the isolation feature in the variance pooling maps.  In such consideration, we adapt the variance statistics for feature classification rather than average correspondence.  
Two distributions of average pooling maps are consistent in the different regions as the smooth patch has smaller average ranges in the TV-regularized reconstruction (Fig.~\ref{fig:synthetic_1D_pooling}(c)) and an oscillatory patch corresponds to the Tikhonov-regularized reconstruction (Fig.~\ref{fig:synthetic_1D_pooling}(d)).  As adaptively referring to two average pooling maps,  we reduce the uncertainty of the exponent distribution from the choice of reconstruction samples.


The exponent is assigned as $p=1$ on the discontinuity regions,  $p=2$ on the oscillatory regions, and the nontrivial exponents ($p \in (1,2$)) are distributed on the smooth regions based on the estimated gradient information from reconstruction samples (Fig.~\ref{fig:example_1_classification_exponent}(b)).  Note that it is not uniform exponent distribution on the flat (i.e., gradient-free) regions on the sparse gradient domain ($[-1,0]$) by reflecting the behavior of reconstruction samples.  Also,  the two 
\begin{figure}[h!]
\centering
\includegraphics[width=1.0\textwidth]{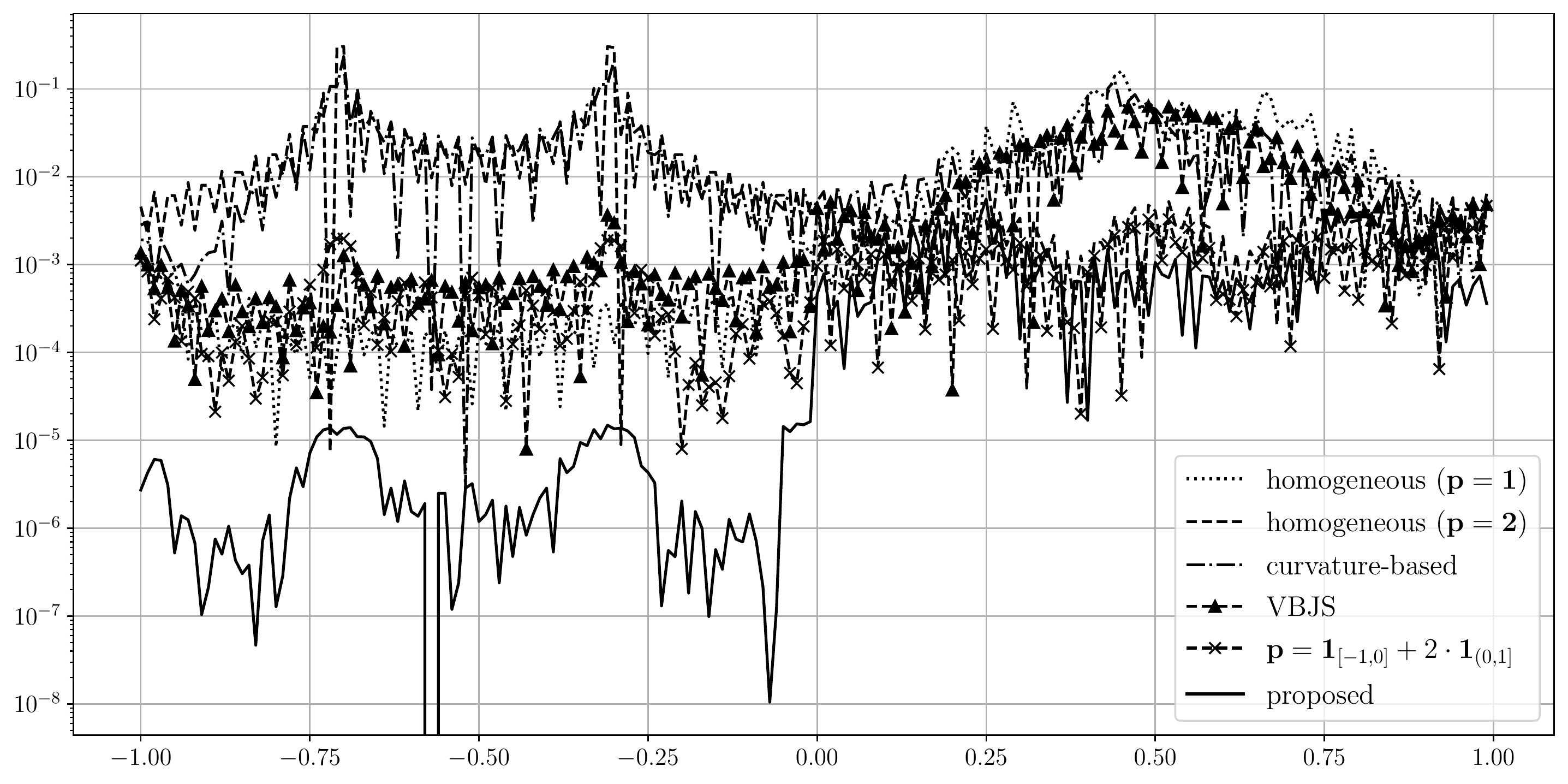}
\caption{\rev{Pointwise errors of six signal recoveries presented in Fig.~\ref{fig:example_1_approximation}}}
\label{fig:example_1_poinwise_errors}
\end{figure}
tails of the oscillation are not classified as oscillations but have relatively large exponents  proportional to the gradient values among the smooth regions.

The signals recovered from different regularizations are presented in Fig.~\ref{fig:example_1_approximation} alongside their point-wise errors in Fig.~\ref{fig:example_1_poinwise_errors}. The recovery from the classical TV (homogeneous $\bm{p}=\bm{1}$) performs well in capturing the discontinuities but is contaminated on the oscillation by flattening the peaks of oscillation as known as the staircase effect (Fig.~\ref{fig:example_1_approximation}(a)).  The case of homogeneous $\bm{p}=\bm{2}$ is the other way around as the oscillation is clearly recovered, but the discontinuity is deteriorated by Gibbs phenomena (Fig.~\ref{fig:example_1_approximation}(b)).  The curvature-based regularization reduces the Gibbs phenomena near the jump discontinuity but is still limited to recover the oscillation fully (Fig.~\ref{fig:example_1_approximation}(c)).  \rev{ The VBJS properly recovers the sparse-gradient regions, but missed the peaks of highly oscillatory signals (\ref{fig:example_1_approximation}(d)).}

As the interpolation idea of $\ell_1$ and $\ell_2$,  the proposed regularization adapts the favorable exponents responding signal features, and the corresponding recovery is accurate through all the regions \rev{(Fig.~\ref{fig:example_1_approximation}(f))}.
We also conduct the experiment that we assign the homogeneous $\bm{p}=\bm{1}$ on $[-1,0]$ and homogeneous $\bm{p}=\bm{2}$ on $(0,1]$ as the exponents are favorable to each separated signal \rev{(Fig.~\ref{fig:example_1_approximation}(e))}. The recovery is comparable to our reconstruction on the oscillatory region;  on the sparse gradient domain $[-1,0]$,  our nonuniform exponent distribution outperforms the uniform distribution $\bm{p}=\bm{1}$ (see Fig.~\ref{fig:example_1_poinwise_errors}).
As the Fourier measurement is a global operator,  the reconstructions on different regions are interrelated. The behavior of reconstruction samples, though possibly inaccurate,  provides a meaningful guideline to design the exponent distribution.

\subsection{Synthetic 2D image recovery} \label{subsec:syn_2d}
We extend our experiment to a 2D image.  We generate the test image
\begin{equation}
u:[-1,1]^2 \rightarrow \mathbb{R}, \hspace{5mm} u(x,y) =
\begin{cases}
\cos\left(18\pi \sqrt{x^2+y^2}\right) & ~\textrm{if}~ \sqrt{x^2+y^2} \leq \frac{4}{9} \\
-18\left(\sqrt{x^2+y^2} -\frac{4}{9}\right)+1 & ~\textrm{if}~ \frac{4}{9}<\sqrt{x^2+y^2} \leq \frac{5}{9} \\
-1 & ~\textrm{if}~ \frac{5}{9}<\sqrt{x^2+y^2} \leq \frac{13}{18}\\
\cos\left(\frac{36}{13}\pi \sqrt{x^2+y^2}\right) & ~\textrm{otherwise}
\end{cases}
\end{equation}
where the features of oscillation,  jump discontinuity,  and ramp are observed.  We consider the discretization with  $128\times 128$ uniform grid on $[-1,1]^2$ (i.e., $N=128^2$). We assume that $25\%$ of original Fourier coefficients are given as the measurement.  The coefficients are sampled from $x$-direction,  and,  as in the previous test,  we choose the first $25\%$ coefficients in increasing order of wavelength.  Also,  the measurement is contaminated by Gaussian noise with $\sigma=18.4$ corresponding to the signal-to-noise ratio (SNR) $25$.

\begin{figure}[!]
\centering
\includegraphics[width=1.0\textwidth]{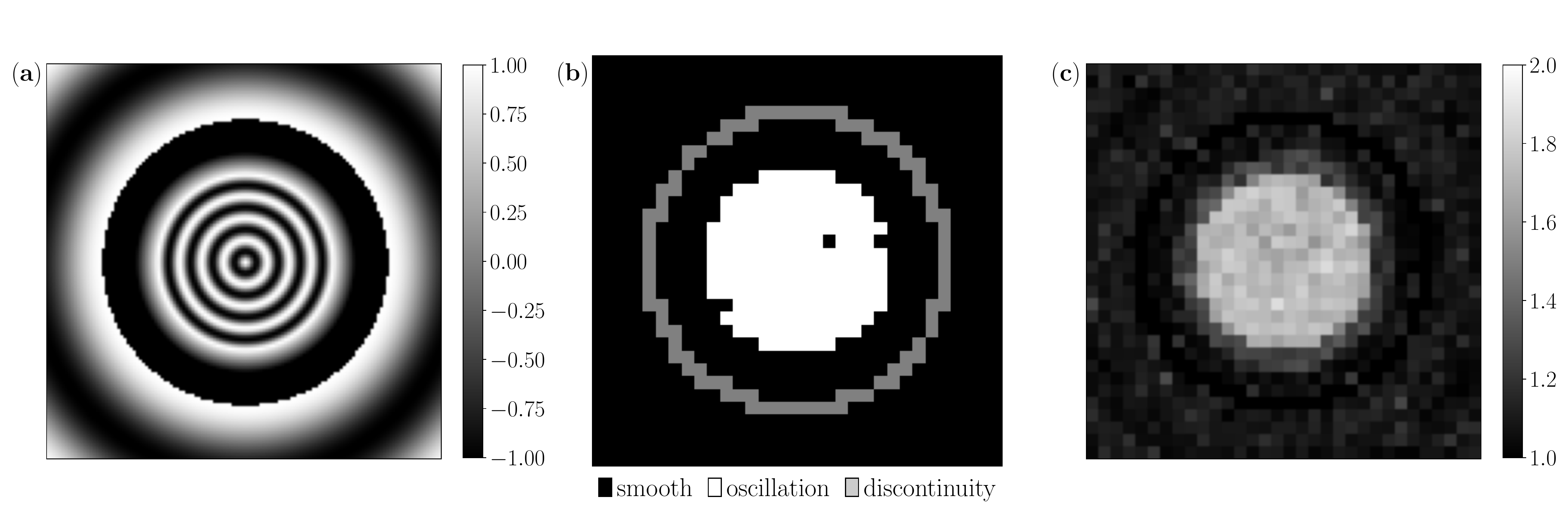}
\caption{True image (a),  local feature classification (b),  exponent distribution (c) with the hyperparameters  $(K,\epsilon, m, c)=(4,1\text{E-}2,2,3.5)$ in algorithm~\ref{algo:wholeProcess}}
\label{fig:example_2_classification_exponent}
\end{figure}

Fig.~\ref{fig:example_2_classification_exponent} presents the feature classification (Fig.~\ref{fig:example_2_classification_exponent}(b)) and the exponent distribution (Fig.~\ref{fig:example_2_classification_exponent}(c)) from algorithms~\ref{algo:patchClassification} and \ref{algo:wholeProcess}.  The smooth region is identified by the collection of patches with low variance in gradient.  Among the patches with high gradient variance,  the discontinuity class is distinguished from the oscillation by observing the directional neighborhoods with smooth patches.  Even under the noisy measurement,  the variance statistics are consistent with identifying the features of regions. In reference to the classification, the exponents $p=1$ are distributed along the discontinuity curve with the patch-sized neighborhood,  and the circular-shaped oscillatory region is covered by the exponents closed to $2$ regarding high gradient distribution in there.  In other ramp regions,  the exponents are adaptively assigned proportionally to the estimated gradient values,  which are close to $1$ in this example.

\begin{figure}[!]
\centering
\includegraphics[width=1.0\textwidth]{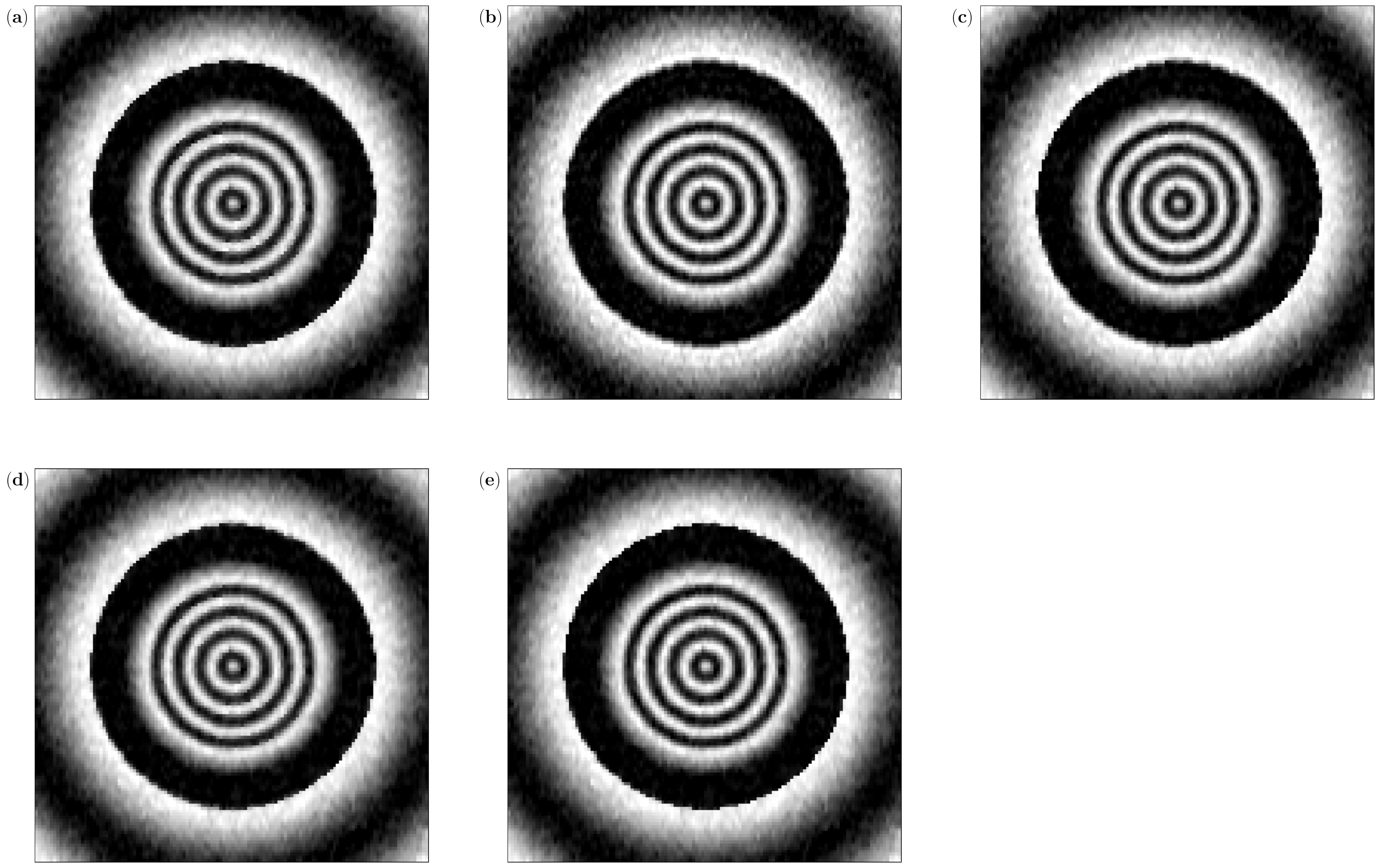}
\caption{\rev{Synthetic 2D image recovery from five different regularizations; (a) homogeneous $\bm{p}=\bm{1}$, (b) homogeneous $\bm{p}=\bm{2}$,  (c) curvature-based regularization,  (d) VBJS regularization,  (e) proposed regularization}}
\label{fig:example_2_approximation}
\end{figure}

\begin{figure}[!]
\centering
\includegraphics[width=1.0\textwidth]{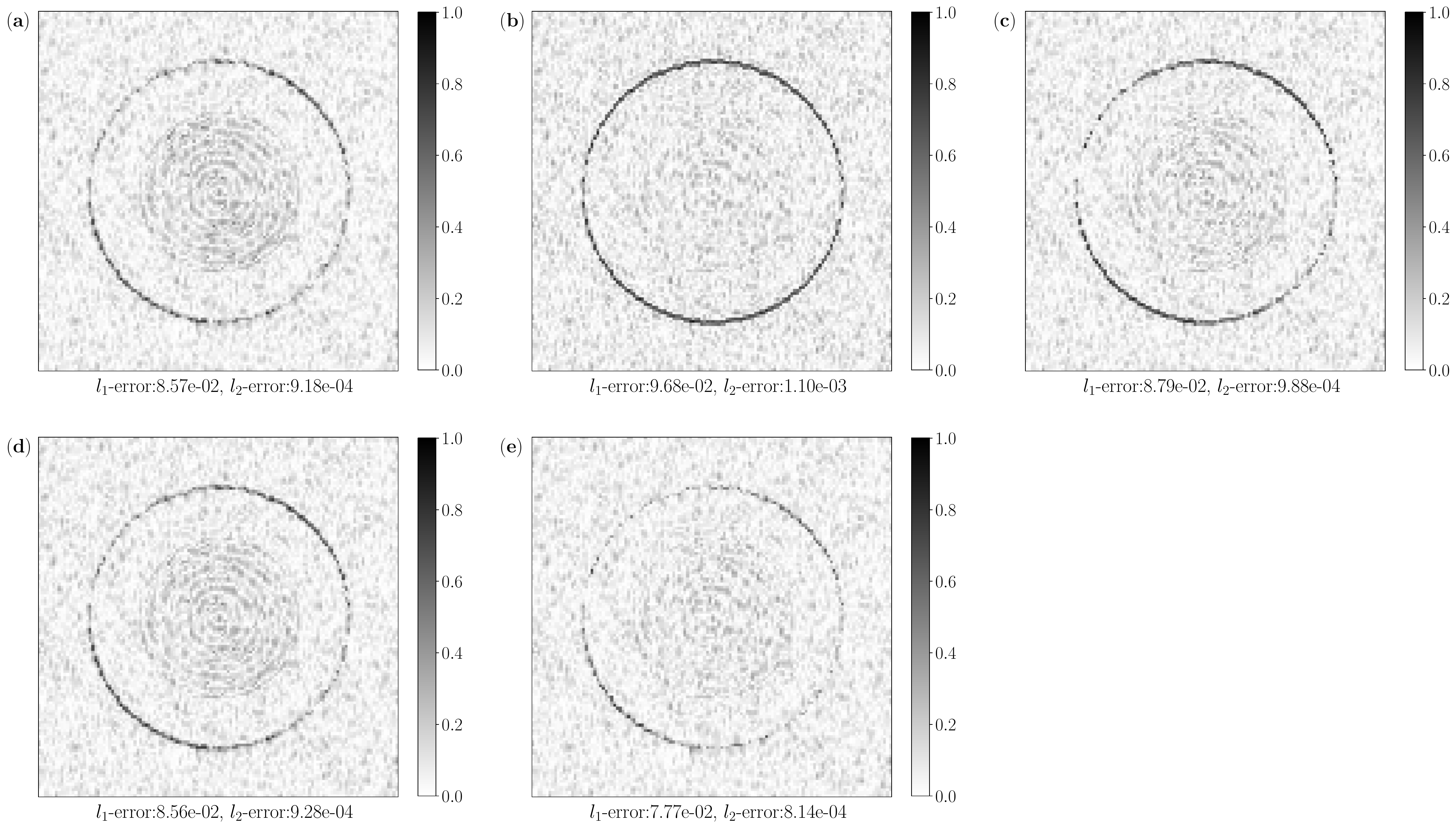}
\caption{\rev{Pointwise errors of image recovery from five different regularizations; (a) homogeneous $\bm{p}=\bm{1}$,  (b) homogeneous $\bm{p}=\bm{2}$, (c) curvature-based regularization,  (d) VBJS regularization, (e) proposed regularization.  The darker color indicates less accurate approximation.}}
\label{fig:example_2_poinwise_errors}
\end{figure}

\begin{figure}[!]
\centering
\includegraphics[width=1.0\textwidth]{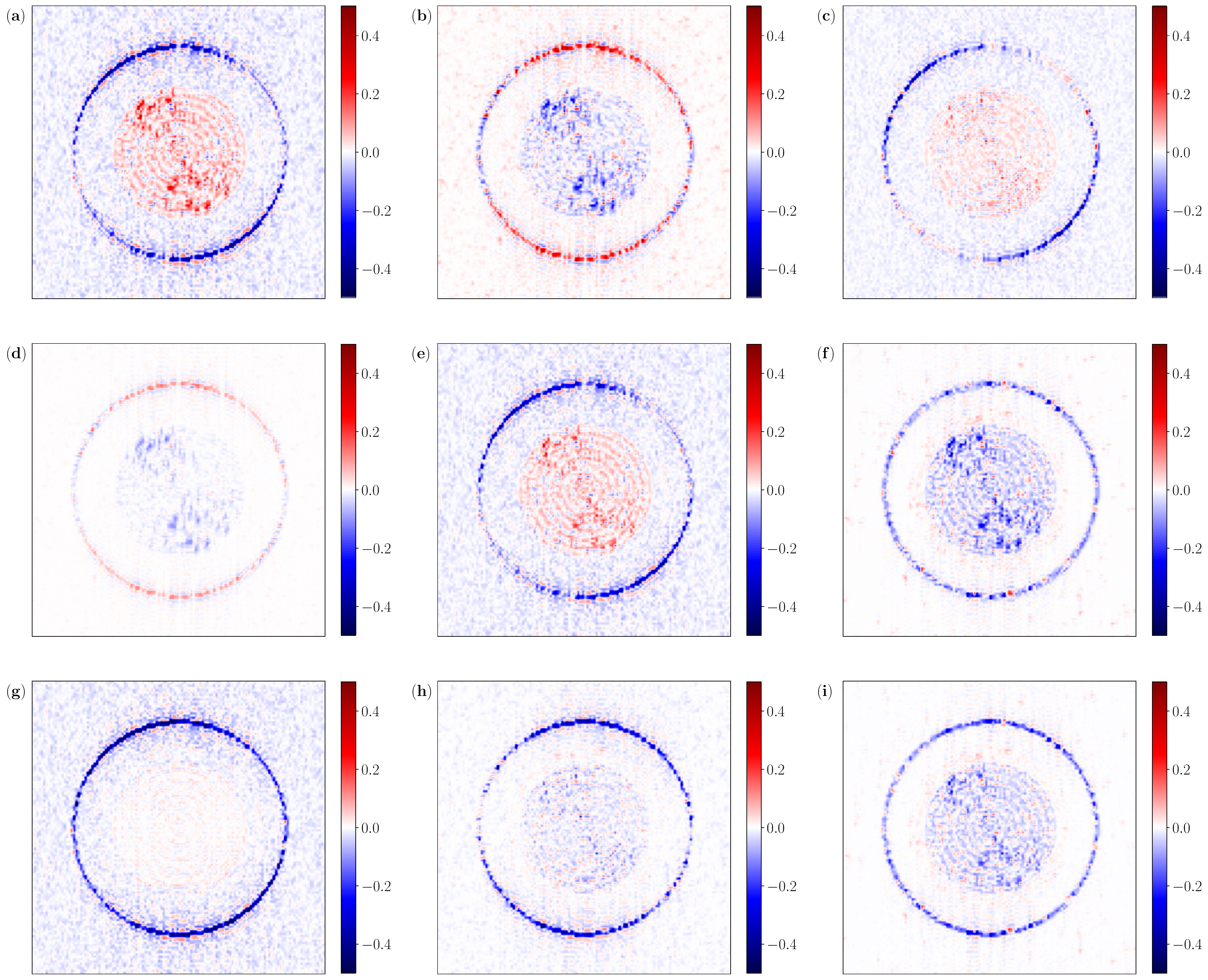}
\caption{\rev{Pointwise error comparison; each plot presents the difference between pointwise errors of two recoveries (a) homogeneous $\bm{p}=\bm{1}$ subtracts homogeneous $\bm{p}=\bm{2}$,  (b) curvature-based regularization subtracts homogeneous $\bm{p}=\bm{1}$,  (c) curvature-based regularization subtracts homogeneous $\bm{p}=\bm{2}$, (d) VBJS regularization subtracts homogeneous $\bm{p}=\bm{1}$,  (e) VBJS regularization subtracts homogeneous $\bm{p}=\bm{2}$,  (f) proposed regularization subtracts homogeneous $\bm{p}=\bm{1}$,   (g)  proposed regularization subtracts homogeneous $\bm{p}=\bm{2}$,  (h) proposed regularization subtracts curvature-based regularization, (i) proposed regularization subtracts VBJS regularization.  The color blue indicates the minuend has less pointwise error than the subtrahend,  and the color red indicates contrariwise.}}
\label{fig:example_2_poinwise_errors_comparison}
\end{figure}

The recovery approximations from the different regularizations are presented in Fig.~\ref{fig:example_2_approximation}, and each pointwise error and pairwise comparisons are summarized in Fig.~\ref{fig:example_2_poinwise_errors} and Fig.~\ref{fig:example_2_poinwise_errors_comparison}.  The homogeneous $\bm{p}=\bm{1}$ (classical TV) is well-performed in the recovery of jump discontinuity, 
and the homogeneous $\bm{p}=\bm{2}$ is corresponding to the oscillatory region.  Fig.~\ref{fig:example_2_poinwise_errors_comparison}(a) presents the comparison of two regularizations, and the advantage of each regularization is distinguishable.   \rev{The curvature-based regularization and VBJS have the trade-off performance on the oscillatory region and other regions including discontinuity and ramps in comparison to homogeneous regularizations (Fig.~\ref{fig:example_2_poinwise_errors_comparison}(b)-(e)).}
The proposed regularization accepts the benefits of two homogeneous regularizations selectively as it is superior to homogeneous $\bm{p}=\bm{1}$ on the oscillatory region \rev{(Fig.~\ref{fig:example_2_poinwise_errors_comparison}(f))} and homogeneous $\bm{p}=\bm{2}$ on the discontinuity \rev{(Fig.~\ref{fig:example_2_poinwise_errors_comparison}(g))}.  On the ramp regions,  the recovery is comparable to the homogeneous $\bm{p}=\bm{1}$ as the exponent distribution is close to $1$ \rev{(Fig.~\ref{fig:example_2_poinwise_errors_comparison}(f))}.  \rev{Compared to the curvature-based regularization and VBJS,  all three regularizations are comparable on the ramp regions,  but the proposed regularization has better performance on the regions of  the discontinuity and the oscillation \rev{(Fig.~\ref{fig:example_2_poinwise_errors_comparison}(h),(i))}. }Overall,  as the signal has oscillatory features in $15.8\%$,   discontinuity in $8.6\%$, and the ramp in $75.6\%$ over the domain,  the proposed regularization improves the homogeneous $\bm{p}=\bm{1}$ with $9.2\%$ in $\ell_1$-error ($11.4\%$ in $\ell_2$-error), the homogeneous $\bm{p}=\bm{2}$ with  $19.6 \%$  in $\ell_1$-error ($25.9\%$ in $\ell_2$-error) (Fig.~\ref{fig:example_2_poinwise_errors}).

\subsection{Real 2D image recovery}
\begin{figure}[!]
\centering
\includegraphics[width=1.0\textwidth]{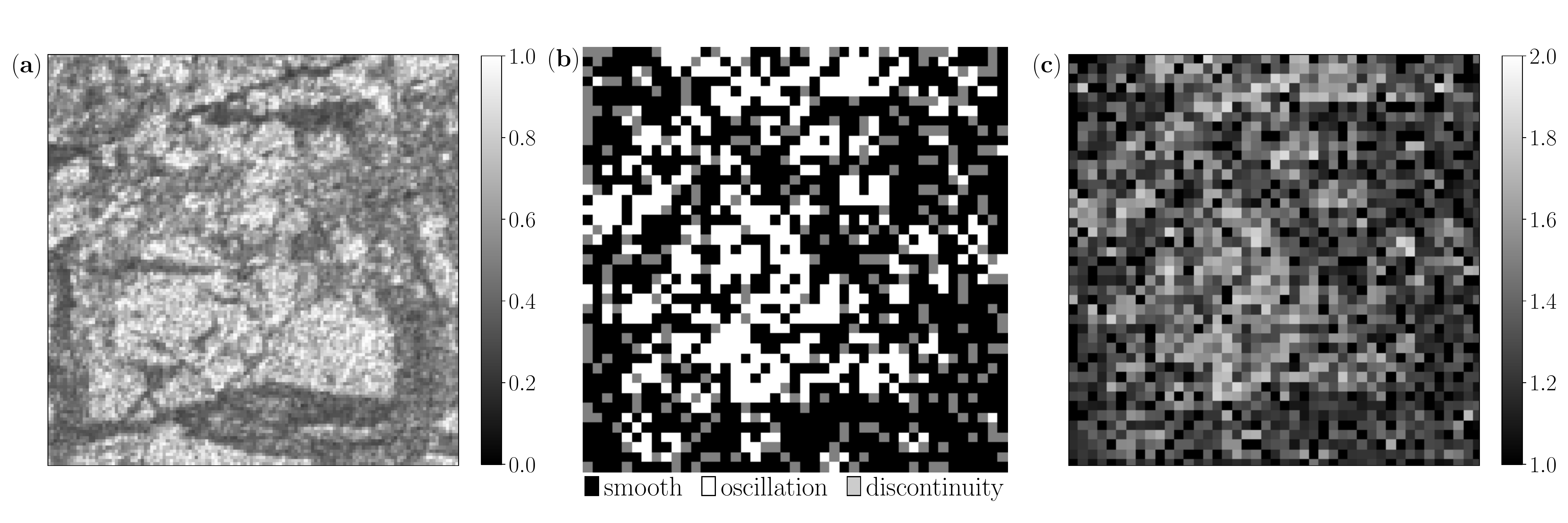}
\caption{True image (a),  local feature classification (b),  exponent distribution (c) with the hyperparameters $(K,\epsilon, m, c)=(3,3\times 10^{-2},2,2)$ in algorithm~\ref{algo:wholeProcess} }
\label{fig:example_3_classification_exponent}
\end{figure}

\begin{figure}[!]
\centering
\includegraphics[width=0.9\textwidth]{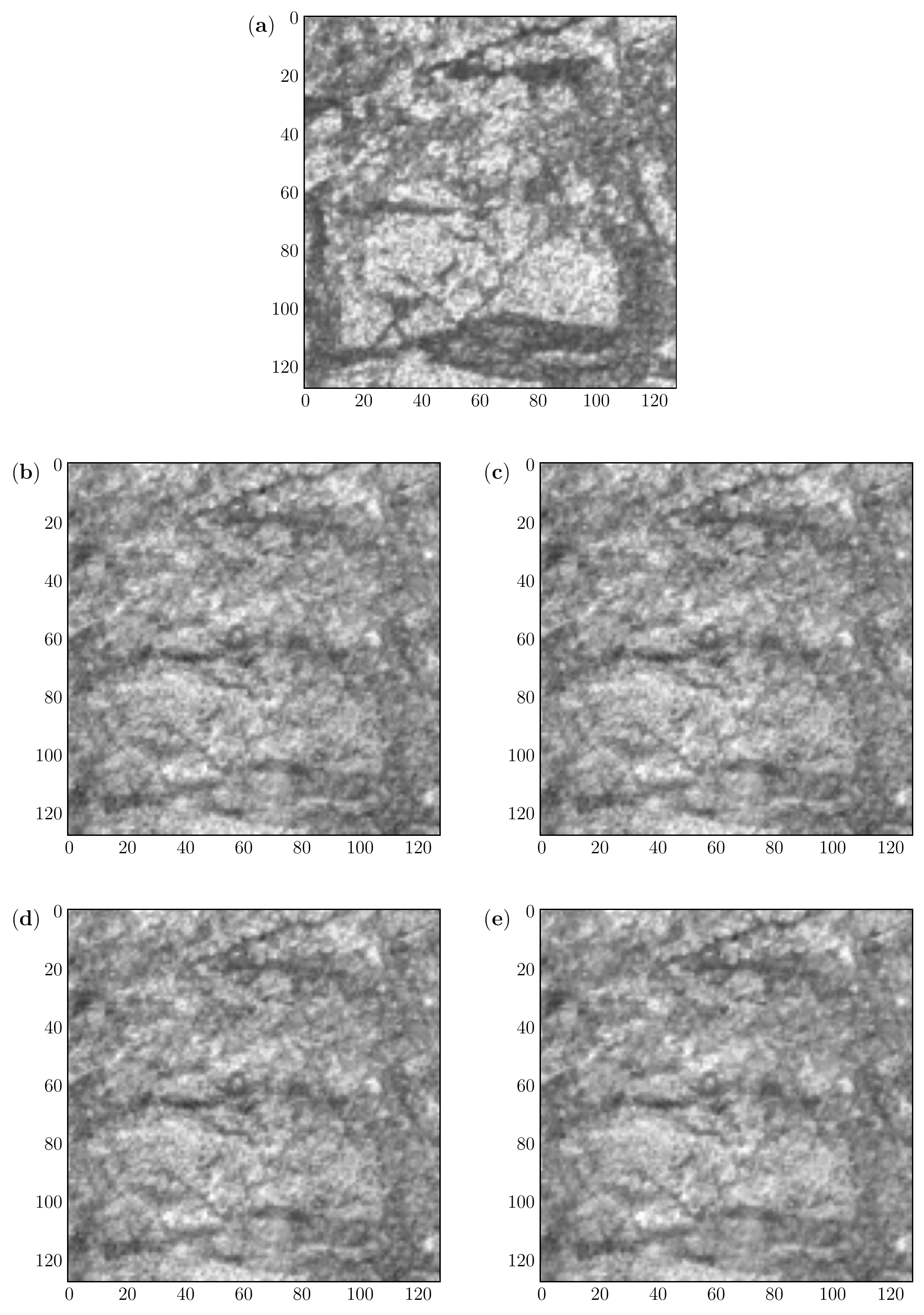}
\caption{Real 2D image recovery from 4 different regularizations; (a) true image,  (b) homogeneous $\bm{p}=\bm{1}$, (c) homogeneous $\bm{p}=\bm{2}$,  (d) curvature-based regularization,  (e) proposed regularization}
\label{fig:example_3_approximation}
\end{figure}

Finally, we test our algorithm for real image reconstruction.  The image is employed from the dataset of Navy Ice Camp See Dragon (ICEX2020).  We downsample the original image and use the part with the resolution $128 \times 128$,  which corresponds to the sea area of size $19.2km\times 19.2km$ (Fig.~\ref{fig:example_3_classification_exponent}(a) or Fig.~\ref{fig:example_3_approximation}(a)). The measurement is the $33\%$ downsampling from the original Fourier coefficients by choosing every three coefficients in $x$-direction.  Also, the measurement is assumed to have additive Gaussian noise with $\sigma=2.08$ corresponding to the signal-to-noise ratio (SNR) 50.

The region classification and the exponent distribution are presented in Fig.~\ref{fig:example_3_classification_exponent}.  The classification map is less readable than the case of the synthetic image, as the features of interest are not simply discernible from the given image.  The classification indicates $15\%$ discontinuity regions and $28\%$ oscillatory regions over the given image.

\begin{figure}[!]
\centering
\includegraphics[width=1.0\textwidth]{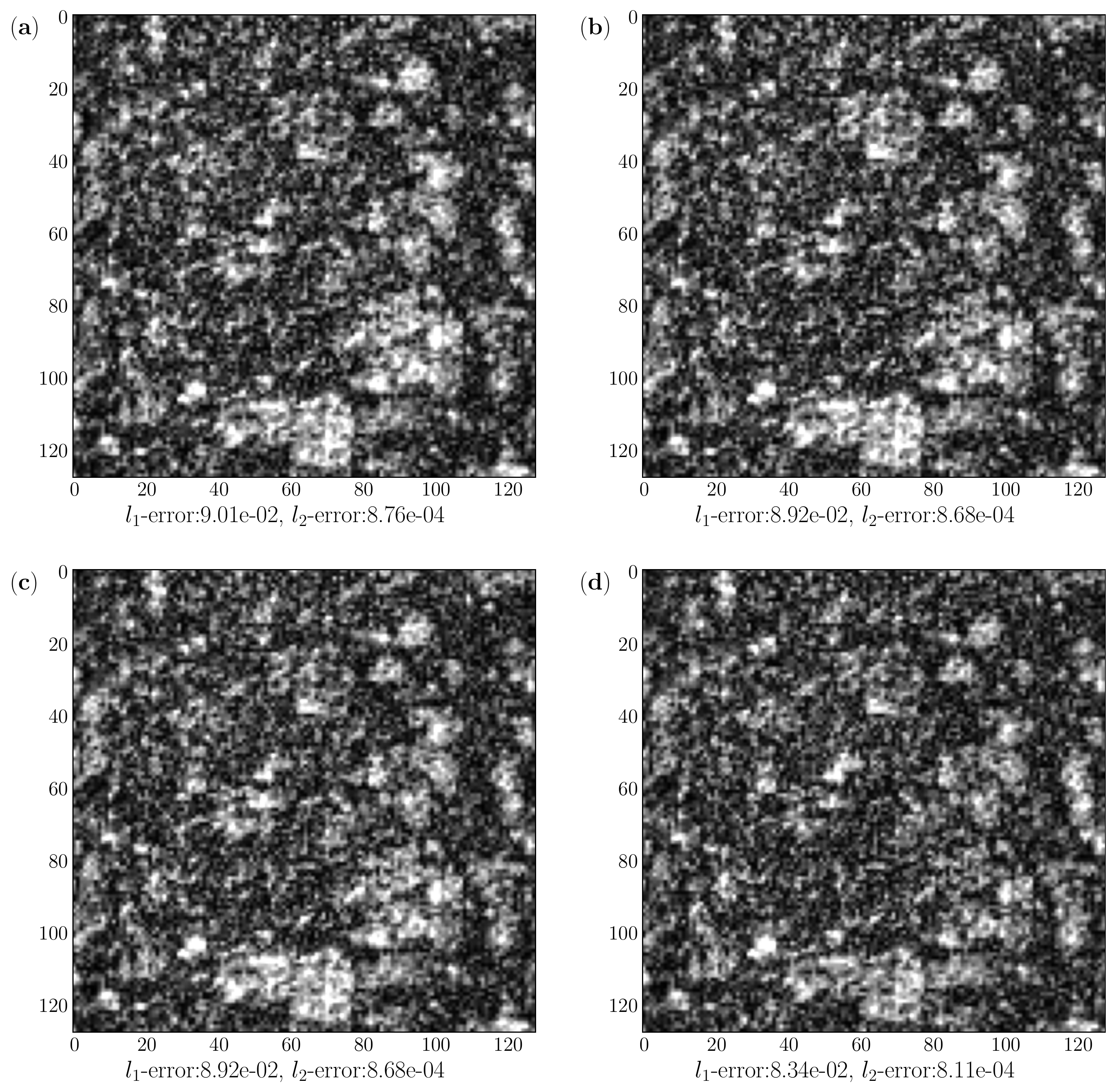}
\caption{Pointwise errors of image recovery from 4 different regularizations; (a) homogeneous $\bm{p}=\bm{1}$,  (b) homogeneous $\bm{p}=\bm{2}$, (c) curvature-based regularization,  (d) proposed regularization}
\label{fig:example_3_poinwise_errors}
\end{figure}

Fig.~\ref{fig:example_3_approximation} displays the reconstruction from different regularizations, and the pointwise errors of each approximation are summarized in Fig.~\ref{fig:example_3_poinwise_errors}.  The recovery from the proposed regularization has resulted in overall improvement as indicated in $\ell_1$ and $\ell_2$ errors.  In particular,  it is better in capturing the underlying features of sea ice textures, which is especially discernible as illustrated by the pointwise error on the region $[80,100] \times [80,100]$ (see Fig.~\ref{fig:example_3_poinwise_errors}).


\section{Conclusions}\label{sec:conclusions}
For an underdetermined inverse problem, regularization plays an important role in stabilizing the inverse problem. 
The current study has focused on the design of the inhomogeneous exponent when the first and the second derivative of the true signal are not available due to indirect and incomplete measurement data. The method proposed in the current study generates multiple reconstructions using standard homogeneous regularization to extract statistical information. The local characteristics of the unknown signal are estimated on a set of patches for the average and the variance pooling of the gradient. The method classifies each patch in an unsupervised way where the exponent value is assigned based on the classification. To efficiently solve the inhomogeneous regularization problem, we has modified the ADMM method that maintains the computational efficiency of the standard ADMM.
A suite of image reconstruction problems in 1D and 2D, including sea ice reconstruction, have been provided to validate the robustness and effectiveness of the proposed method.
The numerical results show the robust performance of the proposed method in comparison with other regularization methods.

In the current study, we focused on incomplete Fourier measurements, which is a global operator. Although each patch is classified independently, signal components in different patches are connected through the measurement operation. We expect that this implicit connection across different patches will affect the patch-wise exponent estimation. We plan to investigate the effect of varying patch sizes and shapes to handle the interconnection across patches and other measurement types, which we leave as future work.

In our numerical tests, the inhomogeneous regularization shows robust performance without weighting. However, it is natural to investigate the interplay between the weighted and the inhomogeneous regularizations and develop a strategy to design the weights. The inhomogeneous regularization changes the shape of the regularization constraint. On the other hand, the weighted regularization changes the length of a simplex while maintaining the geometry of the simplex. As these two approaches have different geometrical interpretations, we expect that each approach has intrinsic performance difference for a certain class of problems, will be reported in another place.

\section*{Acknowledgments}
The authors thank Chris Polashenski for providing the sea ice data. Yoonsang Lee is supported in part by NSF DMS-1912999 and ONR MURI N00014-20-1-2595. 

\bibliographystyle{elsarticle-num}
\bibliography{references}

\end{document}